\documentclass[10pt]{article}
\usepackage{mathrsfs}
\usepackage{tikz}
\usepackage{latexsym,lineno}
\usepackage{epsfig}
\usepackage{color}
\usepackage{amsmath}\usepackage{fleqn}\usepackage{verbatim}\usepackage{epsf}
\usepackage{amsthm}\usepackage{graphicx, float}\usepackage{graphicx}
\usepackage{amsfonts}\usepackage{amssymb}\usepackage{graphpap}
\usepackage{epic}\usepackage{curves}

\newcommand{\be}{\begin{equation}}
\newcommand{\ee}{\end{equation}}
\newcommand{\benum}{\begin{enumerate}}
\newcommand{\eenum}{\end{enumerate}}
\newcommand{\bit}{\begin{itemize}}
\newcommand{\eit}{\end{itemize}}

\begin{document}
\def\s{\subseteq}
\def\n{\noindent}
\def\se{\setminus}
\def\dia{\diamondsuit}
\def\la{\langle}
\def\ra{\rangle}


\title{The ratio of domination and independent domination numbers on trees}

\footnotetext{The first  author is
partially supported by the Summer Graduate Research Assistantship Program of Graduate School. }
\author{ Shaohui Wang$^{1,2}$\footnote{  Corresponding authors: S. Wang (e-mail: shaohuiwang@yahoo.com), B. Wei (e-mail:  bwei@olemiss.edu). } , Bing Wei$^1$ \\ 
\small\emph {1. Department of Mathematics, The University of Mississippi,}\\\small\emph {  University, MS 38677, USA}\\
\small\emph {2. Department of Mathematics and Computer Science, Adelphi University,}\\\small\emph { Garden City, NY 11530, USA}}
\date{}
\maketitle

\begin{abstract} Accepted by Congressus Numernatium(2016).\\

Let $\gamma(G)$ and $i(G)$ be the domination number and the independent domination number of $G$, respectively.  In 1977, Hedetniemi and Mitchell began with the comparison of of $i(G)$ and $\gamma(G)$ and recently
Rad and Volkmann  posted a conjecture that $i(G)/ \gamma(G) \leq \Delta(G)/2$, where $\Delta(G)$ is the maximum degree of $G$. In this work, we prove the conjecture for  trees and provide the graph achieved the sharp bound.

\end{abstract}

\section{Introduction}
Throughout this paper $G = (V, E)$ is a   simple undirected graph with vertex set $  V(G)$ and edge set $ E(G)$.    For  $v\in V(G)$,  $N_G(v)=\{ w\in V(G):  vw\in E(G)\}$ is 
the  open neighborhood of $v$ and $N_G[v]=N_G(v) \cup \{v\}$ is 
the  closed neighborhood of $v$ in $G$.  If $N_G(v) = \phi$,   $v$ is called an isolated vertex.    For $S\subseteq V(G)$, $N_G(S)$ is the open neighborhood of $S$,  $N_G[S]=N_G(S)\cup S$ is the  closed neighborhood of $S$ and  $G-S$ is a subgraph induced by $V(G)-S$.    A graph $F$ is   a forest if it has no cycles.  Specially,  $F$ is a tree if it contains only one component. A double star is a tree with exactly two vertices of  degree greater than 1. In paticular, if the two vertices have same degree, then it is called a balanced double star. 
The line graph $L(G)$ of  a connected graph is a graph such that each vertex of $L(G)$ represents an edge of $G$  and two vertices of $L(G)$ are adjacent if and only if their corresponding edges share a common endpoint in $G$.

It is known that a vertex set $D \subset V(G)$  is a  dominating set if every vertex of $V(G)-D$ is adjacent to some vertices of $D$. The minimum cardinality of a dominating set is called the domination number, denoted by $\gamma(G)$.  Similarly,
a vertex set $I \subset V(G)$  is  an independent dominating set if $I$ is both an independent set and a dominating set in $G$, where an independent set  is a set of vertices in a graph such that  no two of which are adjacent. The minimum cardinality of  an independent dominating set is called the  independent domination number, denoted by $ i(G)$. Currently, lots of  work relating domination number and independent domination number have been studied, referred to surveys \cite{2,4}.

 In 1977,  S. Hedetniemi and S. Mitchell \cite {1977} showed that for any tree T,
 $\frac{ i(L(T))} { \gamma (L(T))} =1$, where $ L(T)$ is the line graph of T.
Because any line graph is a $K_{1,3}$-free graph, 
 R. B. Allan and R. Laskar \cite {1978} extended the previous result  in 1978 and obtained that if a graph  does not have an induced subgraph isomorphic to $K_{1,3}$, then $  i(G)/ \gamma(G) = 1$.
Recently, Goddard et al.\cite{3} considered the ratio $i(G)/ \gamma(G) $ for regular graphs and proved that $i(G)/ \gamma(G) \leq  3/2$ for cubic graphs. In 2013, Southey and Henning \cite{6} improved the previous result to $i(G)/ \gamma(G)  \leq 4/3$ for connected cubic graphs except for $K_{3,3}$.  
During the same year, Rad and Volkmann \cite{5}  got an upper bound of $i(G)/\gamma(G)$ for a graph $G$ and prosposed the conjecture.
 \vskip 2mm {\bf Theorem 1}\emph{ (Rad and Volkmann \cite{5}) 
Let $G$ be a graph, then
$$ \frac{i(G)}{\gamma(G)} \leq
 \left\{ \begin{array}{rcl}
   \frac{\Delta(G)}{2},\;\;\;\;\;\;\;\; \;\;\;\;\;\;\;\;\;\;\;\;\;\;\;&& \mbox{if } 3 \leq \Delta(G) \leq 5,
  \\\Delta(G) -3 +\frac{2}{\Delta(G) - 1},&&   \mbox{if } \Delta(G) \geq 6.
\end{array}\right.$$
}
{\bf Conjecture 2} \emph{(Rad and Volkmann \cite{5}) 
Let $G$  be a graph with $\Delta(G) \geq 3$, then $i(G)/\gamma(G) \leq \Delta(G)/2$.
}

In 2014, Furuta et al.\cite{7} showed that $i(G)/ \gamma(G)  \leq \Delta(G) - 2 \sqrt{\Delta(G)} +2$  for  a graph $G$ and gave the graph achieved the new bound.  However,
when $\Delta(G)$ is big enough, then $ \Delta(G) - 2 \sqrt{\Delta(G)} +2 > \Delta(G) / 2.$ Now there is a natural question that\\
{\it Q:  Is there  other class of graphs,  which has an affirmative answer for Conjecture 2? }

Motivated by Conjecture 2 and the above question, we prove that Conjeture 2 is true for the tree and provide the graph $G $, which  attains the sharp bound $\Delta(G)/2$.
\begin{figure}[thb]
\center
\begin{tikzpicture}
\tikzstyle {every node} =[fill=black!60,circle,inner sep=0.5pt,text=white] \path 
      (2,0) node(0) {$w_1$} 
      (4,0) node(1) {$w_2$}
      (1,1.5) node(2) {$v_1$}
      (1,0.5) node(3) {$v_2$}
      (1, -0.5)node(4){$*$}
      (1,-1.5) node (5) {$v_s$}
      (5,1.5) node(6)  {$u_1$}
      (5,0.5) node (7) {$u_2$}
      (5,-0.5) node (8) {$*$}  
      (5,-1.5) node (9) {$u_s$}   
       ;
\draw[black,line width=1pt]  (0)--(1) (0)--(2) (0)--(3)  (0)--(4) (0)--(5)  
                                       (1)--(6) (1)--(7) (1)--(8) (1)--(9)      
;
\end{tikzpicture}
\caption{A balanced double star}
\label{fig:cor}
\end{figure}
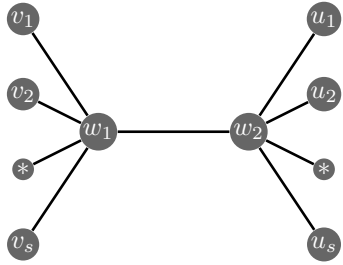

 \vskip 2mm {\bf Theorem 3} \emph{
Let $G$ be a forest, then $$ \frac{i(G)}{\gamma(G)} \leq
 \left\{ \begin{array}{rcl}
1,\;\;\; && \mbox{if } \Delta(G) \leq 2,
  \\\frac{\Delta(G)}{2},&&   \mbox{if } \Delta(G) \geq 3,
\end{array}\right.$$ and the equalities hold if either $\Delta(G) \leq 2$ or each component of $G$ is a balanced double star(see figure 1).
}

As an immediate consequence of Theorem 3, we obtain that 

 \vskip 2mm {\bf Theorem 4} \emph{
Let $G$ be a tree, then $$ \frac{i(G)}{\gamma(G)} \leq
 \left\{ \begin{array}{rcl}
1,\;\;\; && \mbox{if } \Delta(G) \leq 2,
  \\\frac{\Delta(G)}{2},&&   \mbox{if } \Delta(G) \geq 3,
\end{array}\right.$$ and the equalities hold if either $\Delta(G) \leq 2$ or $G$ is a balanced double star(see figure 1).
}

\section{Proof of Theorem 3}

In this section, we will prove Theorem 3 and start with an interesting lemma.

 \vskip 2mm {\bf Lemma 1} \emph{
Let $r_1, r_2,r_3,r_4,t$ be positive numbers with $\frac{r_1}{r_2} \leq t$ and $\frac{r_3}{r_4} \leq t$. Then $\frac{r_1+r_3}{r_2+r_4} \leq t$.  
}

Since $r_1 \leq r_2 t, r_3 \leq r_4 t$, we  replace $r_1$, $r_3$ and obtain that Lemma 1 is true. Next we will give the main proof of this note.

{\bf Proof of Theorem 3.}
For $\Delta(G) \leq 1$, $G$ contains only isolated vertices or edges and $i(G) = \gamma(G)$, that is,  $i(G)/ \gamma(G)  =1$.  Next, we will consider the case of $\Delta(G) \geq 2$ and
  begin with the case that  the forest $G$ contains only one component, that is,  $G$ is a tree.

Let $D$ be a minimum dominating set of $G$. Then $G[D]$ is also 
a forest. We build $\{G_i\}, \{x_i\}$ with $i \geq 1$ as follows:  Let $G_1 = G[D]$ and $ x_1 \in V(G_1)$ with $d_{G_1}(x_1) = 0$ or $1$; For $i \geq 2$, if $V(G_i-N_{G_{i}}[x_{i}]) = \phi$, then stop and set $i = k$. Otherwise, let $G_{i} = G_{i-1}-N_{G_{i-1}}[x_{i-1}]$ and $x_{i} \in V(G_i)$ with $d_{G_i}(x_i) =0 $ or $ 1$.  

Set $X = \{x_1,x_2,...,x_k\}$. Then $X$ is an independent dominating set of $G[D]$
and  $\{N_{G_i}[x_i], 1 \leq i \leq k\}$ is a   partition of $D$, that is, $\sum_{1 \leq i \leq k} (d_{G_i}(x_i) +1) = |D| = \gamma(G)$. Choose $I \subset V(G) - D$ such that 
  $X \cup I$ is an independent dominating set of $G$, that is,  
$
i(G) \leq |X|+|I| = k+ |I|.
$
Since $D$ is a dominating set of $G$, then $I = \cup_{v \in {D-X}}(N_G(v) \cap I) = \cup_{1 \leq i \leq k} (\cup_{v \in {N_{G_i}(x_i)}}(N_G(v) \cap I))$.
By the choice of $x_i$, for $1 \leq i \leq k$ and $v \in N_{G_i}(x_i)$, we have $d_{G_i}(x_i) \leq d_{G_i}(v)$. Thus,
$|N_G(v) \cap I| \leq d_G(v) - d_{G_i}(v) \leq \Delta(G) - d_{G_i}(x_i)$ and
\begin{eqnarray}
|I|   & \leq &\sum_{1 \leq i \leq k}( \sum _{v \in N_{G_i}(x_i)} |N_{G}(v) \cap I|)  \nonumber \\
&\leq  &  \sum_{1 \leq i \leq k} (\sum_{v \in N_{G_i}(x_i)} (\Delta(G) - d_{G_i}(x_i)))  \nonumber \\
&= & \sum_{1 \leq i \leq k} (\sum_{v\in N_{G_i}(x_i)} \Delta (G)) -  \sum_{1 \leq i \leq k} ( \sum_{v \in N_{G_i}(x_i)} d_{G_i}(x_i))  \nonumber\\
&=& (|D| - k ) \Delta(G) - \sum_{1 \leq i \leq k}d_{G_i}(x_i)^2.   
\end{eqnarray}

By $(1)$  and $|D| = \gamma(G)$, we can obtain that  \begin{eqnarray}
i(G) &\leq & k+ |I| \nonumber \\
&\leq &  k+ (|D| - k ) \Delta(G) - \sum_{1 \leq i \leq k}d_{G_i}(x_i)^2  \nonumber \\
&=& \Delta(G)\gamma(G) - \sum_{1 \leq i \leq k } (\Delta(G) -1 +d_{G_i}(x_i)^2), \nonumber
\end{eqnarray}
that is, $$
\frac{i(G)}{\gamma(G)} \leq  \Delta(G) - \frac{ \sum_{1 \leq i \leq k } (\Delta(G) -1 +d_{G_i}(x_i)^2)}{\gamma(G)}.  
$$

Now, it suffices to show that $ - \frac{ \sum_{1 \leq i \leq k } (\Delta(G) -1 +d_{G_i}(x_i)^2)}{\gamma(G)} \leq -\frac{\Delta(G)}{2}$, that is,
\begin{eqnarray} 
\sum_{1 \leq i \leq k} (\Delta(G) - 1 + d_{G_i}(x_i)^2)  &\geq& \frac{1}{2} \Delta(G) \gamma(G)  \nonumber \\&=& \frac{1}{2} \Delta(G) (\sum_{1 \leq i \leq k} (d_{G_i}(x_i) +1)) 
\end{eqnarray}
By the construction of $G_i, x_i$,   $d_{G_i}(x_i) = d_{G_i}(x_i)^2 = 0$ or $1$.
Thus, $(2)$ is the same as $(3)$ below.
\begin{eqnarray} 
 &&\Leftrightarrow
k\Delta(G) - k + \sum_{1\leq i \leq k} d_{G_i}(x_i)  - \frac{1}{2} \Delta(G)(\sum_{1\leq i \leq k} d_{G_i}(x_i) )   \nonumber \\&& \;\;\;\;\; -\frac{1}{2} \Delta(G) k  \geq 0  \nonumber
\\&&
 \Leftrightarrow
(1- \frac{1}{2} \Delta(G)) ((\sum_{1\leq i \leq k} d_{G_i}(x_i) ) - k) \geq 0  
\end{eqnarray}

Furthermore, $d_{G_i}(x_i) = 0$ or $1$ yields that   $(\sum_{1\leq i \leq k} d_{G_i}(x_i) ) - k \leq 0$.
Since $\Delta(G) \geq 2$,  then
 $1- \frac{1}{2} \Delta(G) \leq 0$. Thus, $(3)$ is true, that is, Theorem 3 is true for the tree.

Next we will consider the case that $G$ has more than one component. In this case,  each component of $G$ is either  an isolated vertex or a tree, say $G_1, G_2, ..., G_s$ with an integer $s \geq 2$. For $1 \leq j \leq s $, if $G_j$ is an isolated vertex, then $i(G_j)/\gamma(G_j) = 1/1 \leq \Delta(G)/2$; If $G_j$ is a tree, by the above proof,   $i(G_j)/\gamma(G_j) \leq   \Delta(G)/2$.  Finally,   using Lemma 1,    $i(G)/\gamma(G) \leq  \Delta(G)/2$ holds for the forest. 
 Furthermore, if $\Delta(G) \leq 2$, all  forests achieve the bound; if $\Delta(G) \geq 3$, the union of balanced double stars attain the bound. Thus, Theorem 3 is true.

 $\hfill\Box$

\end{document}